\documentclass{amsart}
\usepackage{amssymb}
\usepackage{color}

\newtheorem{theorem}{Theorem}[section]
\newtheorem{lemma}[theorem]{Lemma}

\theoremstyle{definition}
\newtheorem{definition}[theorem]{Definition}

\theoremstyle{remark}
\newtheorem{remark}[theorem]{Remark}

\numberwithin{equation}{section}

\begin{document}

\title[A Bernstein-type inequality for rational functions]
{A Bernstein-type inequality for rational functions in weighted
Bergman spaces}

\author{Anton Baranov}
\address{Department of Mathematics and Mechanics, Saint Petersburg State University, 28,
Universitetski pr., St. Petersburg, 198504, Russia}

\email{anton.d.baranov@gmail.com}

\author{Rachid Zarouf}

\address{CMI-LATP, UMR 6632, Universit\'e de Provence, 39, rue F.-Joliot-Curie,
13453 Marseille cedex 13, France}
\email{rzarouf@cmi.univ-mrs.fr}

\thanks{The first author is supported by the Chebyshev Laboratory
(St.Petersburg State University) under RF Government grant 11.G34.31.0026
and by RFBR grant 12-01-00434.}



\begin{abstract}Given $n\geq1$ and $r\in[0,\,1),$
we consider the set $\mathcal{R}_{n,\, r}$
of rational functions having at most $n$ poles
all outside of $\frac{1}{r}\mathbb{D},$
were $\mathbb{D}$ is the unit disc of the complex plane. We give
an asymptotically sharp Bernstein-type inequality for functions
in $\mathcal{R}_{n,\, r}\:$
in weighted Bergman spaces with
``polynomially'' decreasing weights.
We also prove that this result
can not be extended to weighted Bergman spaces  with {}``super-polynomially''
decreasing weights.
\end{abstract}

\maketitle

\section{Introduction}

Estimates of the norms of derivatives for polynomials and rational
functions (in different functional spaces) is a classical topic of
complex analysis (see surveys  by A.A. Gonchar \cite{Go},
V.N. Rusak \cite{Ru}, and P. Borwein and T. Erd\'elyi
\cite[Chapter 7]{BoEr}). Such inequalities
have applications in many domains of analysis;
to mention just some of them:
1) matrix analysis and
in operator theory (see {}``Kreiss Matrix Theorem'' \cite{LeTr,Sp}
or \cite{Z1,Z5} for resolvent estimates of power bounded matrices),
2) inverse theorems of rational approximation (see \cite{Da,Pel,Pek}),
3) effective Nevanlinna--Pick interpolation problems (see
\cite{Z3,Z4}).

Here, we present Bernstein-type
inequalities for rational functions $f$ of degree $n$ with
poles in $\left\{ z\,:\;\left|z\right|>1\right\} $, involving Hardy
norms and weighted Bergman norms.
Let $\mathcal{P}_{n}$ be the complex space of polynomials
of degree less or equal to $n\ge1$.
Let $\mathbb{D}=\left\{ z\in\mathbb{C}\,:\left|z\right|<1\right\} $
be the unit disc of the complex plane and $\overline{\mathbb{D}}=\left\{ z\in\mathbb{C}\,:\left|z\right|\leq1\right\} $
its closure. Given $r\in[0,\,1),$ we define
\[
\mathcal{R}_{n,\, r}=\left\{ \frac{p}{q}\,:\; p,\,
q\in\mathcal{P}_{n},\;{\rm d}^{\circ}p<{\rm d}^{\circ}q,\;
q(\zeta) \ne 0\;\; |\zeta| < \frac{1}{r}\right\},
\]
(where ${\rm d}^{\circ}p$ denotes the degree of $p\in\mathcal{P}_{n}$),
the set of all rational functions in $\mathbb{D}$ of degree less
or equal than $n\ge1$, having at most $n$ poles all outside of $\frac{1}{r}\mathbb{D}.$
Notice that for $r=0$, we get $\mathcal{R}_{n,\,0}=\mathcal{P}_{n-1}$.


\subsection{Definitions of Hardy spaces and radial weighted Bergman spaces}
We denote by ${\rm Hol}\left(\mathbb{D}\right)$ the space of all holomorphic
functions on $\mathbb{D}$. From now on, if $f\in{\rm Hol}\left(\mathbb{D}\right)$
then for every $\rho\in(0,\,1)$ we define
\[
f_{\rho}\::\:\xi\mapsto f\left(\rho\xi\right),\quad \xi
\in\frac{1}{\rho}\mathbb{D}.
\]
We consider the two following scales of Banach spaces $X\subset{\rm Hol}\left(\mathbb{D}\right):$

\textbf{a.} The Hardy spaces $H^{p}=H^{p}(\mathbb{D}),$
$1\leq p\leq\infty:$
\[
H^{p}=\left\{ f\in{\rm Hol}\left(\mathbb{D}\right):
\:\left\Vert f\right\Vert _{H^{p}}^{p}=\sup_{0\leq\rho<1}
\int_{\mathbb{T}}\left|f_{\rho}(\xi)\right|^{p}{\rm d}m(\xi)<\infty\right\} ,
\]
where $m$ stands for the normalized Lebesgue measure on $\mathbb{T}=
\left\{ z\in\mathbb{C}:\:\left|z\right|=1\right\}$. As usual,we denote by
$H^\infty$ the space of all bounded analytic functions in $\mathbb{D}$.

\textbf{b.} The radial weighted Bergman spaces $L_{a}^{p}\left(w\right)$,
$1\leq p<\infty$ (where $"a"$ means analytic),
\[
L_{a}^{p}\left(w\right)=\left\{ f\in{\rm Hol}\left(\mathbb{D}\right):\:\left\Vert f\right\Vert _{L_{a}^{p}\left(w\right)}^{p}=\int_{0}^{1}\rho w\left(\rho\right)\int_{\mathbb{T}}\left|f_{\rho}\left(\zeta\right)\right|^{p}{\rm d}m(\zeta){\rm d}\rho<\infty\right\} ,
\]
where the weight $w$ satisfies $w\geq0$ and $\int_{0}^{1}
w(\rho) {\rm d}\rho<\infty.$
For the classical power weights
$w(\rho)=w_{\beta}(\rho)=\left(1-\rho\right)^{\beta}$,
$\beta>-1$, we have $L_{a}^{p}\left(w_{\beta}\right)=
L_{a}^{p}\big( (1-\left|z\right|)^{\beta}{\rm d}A(z)\big)$,
$A$ being the normalized area measure on $\mathbb{D}$.

For general properties of these spaces we refer to \cite{HedKorZhu,Zhu}.

From now on, for two positive functions $a$ and $b$, we say that
$a$ is dominated by $b$, denoted by $a \lesssim b$, if there is a constant
$c>0$ such that $a\leq cb;$ and we say that $a$ and $b$ are comparable,
denoted by $a\asymp b$, if both $a \lesssim b$ and $b \lesssim a$.


\subsection{Statement of the problem and known results}
By Bernstein-type inequalities for rational functions
one usually understands the inequalities of the form

\begin{equation}
\label{0}
\|f'\|_X \le \phi_{X,\, Y}(n) \|f\|_Y, \qquad f\in \mathcal{R}_n,
\end{equation}

where $\mathcal{R}_n$ is the set of all proper
rational functions of degree at most $n$
with the poles in  $\{|z|>1\}$,
$X$ and $Y$ are some normed spaces of functions analytic in the unit disc,
and $\phi$ is some increasing (often polynomially growing) function.
Thus, for a given pair of the function spaces
$X$ and $Y$, the question is to determine the dependence on $n$
for the norm of the differentiation operator $(\mathcal{R}_n, \|\cdot\|_X)$
to $Y$. Bernstein-type inequalities of E.P. Dolzhenko
\cite{Dol} and A.A. Pekarskii \cite{Pek} are of this form; e.g., it is shown in
\cite{Dol} that
$$
\|f\|_{H_1^1} \le c_1 n \|f\|_\infty,  \qquad
\|f\|_{B_{2,2}^{1/2}} \le    c_2 n ^{1/2} \|f\|_\infty,
\qquad f\in \mathcal{R}_n,
$$
where $H_1^1$ is the Hardy--Sobolev space, and $B_{2,2}^{1/2}$
is the Besov (or Dirichlet) space, see the definition in Section 3.
Let us also mention that this problem is a part of a
more general one given by G. Lorentz in a letter sent to T. Erd\'elyi
in 1988 (see \cite{Erd}).

Looking at (\ref{0}), we notice that for some choices of $X$ and $Y$,
we have $\phi_{X,\, Y}(n)=+\infty$ for every $n=1,\,2,\,\dots$.
Indeed, it may happen for instance when the poles of our function
$f$ are allowed to be arbitrary close to the torus $\mathbb{T}$
: we can observe this phenomenon for example in the special case $X=Y=H^{p},\,1\leq p\leq+\infty$
but also when $X=Y=L_{a}^{p}(w),\,1\leq p\leq+\infty.$ This observation
leads us to come back on the problem in  (\ref{0}) and to state it more generally
: that is replacing $\mathcal{R}_{n}$ by $\mathcal{R}_{n,\, r}$
(for any fixed $r\in[0,\,1)$) and $\phi_{X,\, Y}(n)$ by $\phi_{X,\, Y}(n,\, r)$
so that to focus on this phenomenon of {}``natural dependence on
the parameter $r$''. For most of the classical cases already studied
by others (for instance E. P. Dolzhenko \cite{Dol}, A. A. Pekarskii \cite{Pek}, V.V.
Peller \cite{Pel}) the spaces $X$ and $Y$ are such that $\sup_{r\in(0,\,1)}\phi_{X,\, Y}(n,\, r)<+\infty$:
in this case we can set $\phi_{X,\, Y}(n)=\sup_{r\in(0,\,1)}\phi_{X,\, Y}(n,\, r)$.
As a consequence, if $\sup_{r\in(0,\,1)}\phi_{X,\, Y}(n,\, r)=+\infty,$
it may be of interest to search (as a continuation of the investigations
of the second author \cite{Z2, Z6}) for the {}``best possible'' $\phi_{X,\, Y}(n,\, r)$
in an asymptotically sense, that is to say as $n\rightarrow\infty$
and $r\rightarrow1^{-}$. This question has already been answered
for the case $X=Y=H^{p},\,1\leq p\leq+\infty$ by K. M. Dyakonov \cite{Dy2}
see (\ref{hardy}) below. In this paper, we  answer the same question
for the case $X=Y=L_{a}^{p}(w),\,1\leq p\leq+\infty.$
Let us give a general formulation of our problem for the special case
$X=Y$ for which we set $\mathcal{C}_{n,\, r}(X)=\phi_{X,\, Y}(n,\, r)$
: given a Banach space
$X$ of holomorphic functions in $\mathbb{D},$ we are searching for
the best possible constant $\mathcal{C}_{n,\, r}(X)$ such
that
\[
\left\Vert f'\right\Vert _{X}\leq\mathcal{C}_{n,\, r}(X)
\left\Vert f\right\Vert _{X},
\qquad f\in\mathcal{R}_{n,\, r}.
\]

For the case where $X=H^p$ is a Hardy space, an estimate
which gives a correct order of growth for $\mathcal{C}_{n,\, r}(X)$
was obtained by K.M. Dyakonov \cite{Dy2}
(as a very special case of more general
results): for any $p\in [1, \infty]$ there exist positive constants $A_p$
and $B_p$ such that
\begin{equation}
\label{hardy}
A_p \frac{n}{1-r} \le \mathcal{C}_{n,\, r}(H^p) \le B_p \frac{n}{1-r}
\end{equation}
for all $n\geq1$ and $r\in[0,\,1)$.
More precisely, the upper estimate for $p\in(1,\,+\infty)$ is
treated in \cite[Theorem 1]{Dy2}, the case $p=1,$ in \cite[Corollary 1]{Dy2},
and the case $p=+\infty$ (known much earlier)
is given in \cite[Theorem 7.1.7]{BoEr}.
The below estimate follows trivially when applying the
differentiation operator to  the
test function $f(z) = (1-rz)^{-n}$.

For the case $p=2$ an asymptotically sharp result was obtained
later in \cite{Z2}: for any $r \in (0,1)$ there exists the limit
$$
\lim_{n\rightarrow\infty}\frac{\mathcal{C}_{n,\, r}\left(H^{2}\right)}{n}
=\frac{1+r}{1-r}.
$$

Related results about Bernstein-type inequalities
in a more general setting of the
so-called model or star invariant subspaces
may be found in \cite[Theorems 10,11]{ind}, \cite[Theorem 1]{Dy1},
and \cite{bar1, bar2}.


\subsection{Main results}

We obtain estimates for the derivatives of rational functions
with respect to weighted Bergman norms.
It turns out that there
is an essential difference between slowly (polynomially)
decreasing weights and fast (super-polynomially) decreasing
weights. In the first case we have a two-sided estimate
analogous to (\ref{hardy}), while in the second case only the above
estimate remains true. Let us give the precise definitions.
Recall that $w$ is always an integrable nonnegative function on $(0,\,1)$.

\begin{definition}

(\textit{Polynomially decreasing weights}) The weight $w$ is said
to be \textit{$\gamma$-polynomially decreasing} if there exists $\gamma>0$
such that
\[
\rho\mapsto(1-\rho)^{-\gamma}w(\rho),
\]
is increasing on $\left[r_{0},\,1\right)$ for some $\;0\leq r_{0}<1$.
We say that $w$ is \textit{polynomially decreasing} if
it is \textit{$\gamma$-polynomially
decreasing for some $\gamma>0.$ }

\end{definition}

\begin{definition}
(\textit{Super-polynomially decreasing weights}) The weight $w$ is said
to be \textit{super-polynomially decreasing} if for any $\gamma>0$ there
exists $r(\gamma)\in(0,1)$ such that the function
\[
\rho\mapsto(1-\rho)^{-\gamma}w(\rho),
\]
decreases on the interval $[r(\gamma),\,1)$.

\end{definition}

Typical example of the weights from the first class are given by
$w(r) = (1-r)^{\beta}$, $\beta>-1$,
or  $w(r) = (1-r)^{\beta}(|\log (1-r)|+1)^{\gamma}$,
$\beta>-1$, $\gamma\in \mathbb{R}$. The weights $w(r) =
\exp\big( -c(1-r)^{-\gamma} \big)$, $c>0$, $\gamma>0$
are super-polynomially decreasing.

Our first result may be considered as an analogue of Dyakonov's theorem
for the radial weighted Bergman spaces.

\begin{theorem}
\label{ourth1}
Let $1\leq p<\infty$ and let $w$ be
an integrable nonnegative function on $[0,\,1)$.
Then there exists a positive constant $K$
depending only on $p$ (but not on the weight $w$) such that
\begin{equation}
\label{1}
\mathcal{C}_{n,\, r}\left(L_{a}^{p}\left(w\right)\right)\leq K\frac{n}{1-r}
\end{equation}
for all $r\in[0,\,1)$ and $n\geq1$.
Moreover, if we fix $r \in (0,1)$ and let $n$ tend
to infinity, then we have
\begin{equation}
\label{1-1}
  \frac{\widetilde{K}r}{1-r}\le\liminf_{n\to\infty}
  \frac{C_{n,r}(L_{a}^{p}(w))}{n}\le\limsup_{n\to\infty}\frac{C_{n,r}(L_{a}^{p}(w))}{n}\le\frac{K}{1-r},
\end{equation}
where $\widetilde{K}$ is, as $K$, a positive constant depending
only on $p$.
\end{theorem}

The next theorem shows that for the polynomially decreasing weights
the quantity $\mathcal{C}_{n,\, r} (L^p_a(w))$ admits a below estimate
of the same form.

\begin{theorem}
\label{ourth2}
If $w$ is $\gamma$-polynomially
decreasing\textit{,} then there exists a positive constant $K'$ depending
only on $w$ and $p$ such that
\begin{equation}
\label{2}
K'\frac{n}{1-r}\leq\mathcal{C}_{n,\, r}\left(L_{a}^{p}\left(w\right)\right)\leq K\frac{n}{1-r},
\end{equation}
where $K$ is defined in \eqref{1} and where the left-hand side
inequality of \eqref{2} holds for all $r\in[0,\,1)$ and
$n \ge \frac{\gamma+3}{p}+1$.
In particular, \eqref{2} holds for the classical weights
$w(\rho)=w_{\beta}(\rho)=\left(1-\rho\right)^{\beta}\rho$, $\beta>-1.$
\end{theorem}

The polynomial decrease is essential and provides a sharp
bound for the validity of the uniform estimate (\ref{2}) for all possible values
of $n$ and $r$. Namely, if the weight
is super-polynomially decreasing, then (\ref{2})
will fail along some sequence of radii.

\begin{theorem}
\label{ourth3}
Suppose that $w$ is super-polynomially decreasing.
Then there exists a sequence $r_{n} \to 1-$ such that for any $p$, 
$$
\frac{C_{r_{n},n}(L_{a}^{p}(w))}{n}=o\bigg(\frac{1}{1-r_{n}}\bigg),
\qquad n\to\infty.
$$
\end{theorem}

{\bf Acknowledgements.}
The authors are  deeply grateful
to Nikolai Nikolski, Alexander Borichev, and Evgueny Doubtsov
for many helpful discussions, constructive comments and
precious remarks which definitively helped to improve the manuscript.


\section{Proofs of Theorems \ref{ourth1} and \ref{ourth2}}

\begin{proof} [Proof of Theorem \ref{ourth1}]
First, we notice that for any $0\leq\alpha<1$,

\begin{equation}
  \label{8}
  \left\Vert f\right\Vert _{L_{a}^{p}\left(w\right)}^{p}
  \asymp\iint_{u\in C_{\alpha}}\rho\left|f(\rho\xi) \right|^{p}
  w(\rho) {\rm d}m(\zeta){\rm d}\rho
\end{equation}

for all $f\in L_{a}^{p}(w),$ where $C_{\alpha}=
\left\{ z\,:\;\alpha<\left|z\right|<1\right\}$.
Let $f\in\mathcal{R}_{n,\, r}$ with $r\in[0,\,1)$
and $n\geq1$. Using \eqref{8} with $\alpha=\frac{1}{2}$ we get
\begin{eqnarray*}
\left\Vert f'\right\Vert _{L_{a}^{p}\left(w\right)}^{p} &
\asymp & \iint_{\rho\xi\in C_{1/2}}\left|f'(\rho\xi)\right|^{p}
w\left(\rho\right){\rm d}m(\zeta){\rm d}\rho  \\
 & = & \int_{\frac{1}{2}}^{1}\rho\, w(\rho)
\frac{1}{\rho^{p}}\left(\left\Vert \left(f_{\rho}\right)^{'}
\right\Vert _{H^{p}}^{p}\right){\rm d}\rho.
\end{eqnarray*}

Now using the fact that $f_{\rho}\in\mathcal{R}_{n,\,\rho r}
\subset\mathcal{R}_{n,\, r}$ for every $\rho\in(0,\,1)$, we get
\begin{eqnarray*}
\int_{\frac{1}{2}}^{1}\rho w(\rho)\frac{1}{\rho^{p}}\left(\left\Vert \left(f_{\rho}\right)^{'}\right\Vert _{H^{p}}^{p}\right){\rm d}\rho & \leq & \left(2\mathcal{C}_{n,\, r}\left(H^{p}\right)\right)^{p}\int_{\frac{1}{2}}^{1}\rho w\left(\rho\right)\left(\left\Vert f_{\rho}\right\Vert _{H^{p}}\right)^{p}{\rm d}\rho\\
 & \asymp & \left(\mathcal{C}_{n,\, r}\left(H^{p}\right)\right)^{p}\left\Vert f\right\Vert _{L_{a}^{p}\left(w\right)}^{p}.
\end{eqnarray*}
In particular, using the right-hand side inequality of \eqref{hardy},
we get
\[
\mathcal{C}_{n,\, r}\left(L_{a}^{p}\left(w\right)\right)\leq K_{p}\frac{n}{1-r},
\]
 for all $p\in[1,\,\infty),$ and $\beta\in(-1,\,\infty)$, where
$K_{p}$ is a constant depending on $p$ only.

Now, let us prove \eqref{1-1}. Let
\[
f_{n}(z)=\frac{1}{(1-rz)^{n}}\in\mathcal{R}_{n,r},
\]
and $D=\{z\in\mathbb{D}:\:|1-rz|\le2|1-r|\}$. We claim that
\[
\|f_{n}\|_{L_{a}^{p}(w)}^{p}\sim\int_{D}|f_{n}(z)|^{p}w(z)\mbox{d}A(z),\qquad n\to\infty,
\]
and, analogously,
\[
\|f'_{n}\|_{L_{a}^{p}(w)}^{p}\sim\int_{D}|f_{n}'(z)|^{p}w(z)\mbox{d}A(z),
\qquad n\to\infty.
\]
Indeed, by a very rough estimate
\[
\int_{\mathbb{D}\setminus D}|f_{n}(z)|^{p}w(z)\mbox{d}A(z)\le\frac{C_{1}}{2^{pn}(1-r)^{pn}},
\]
where $C_{1}>0$ depends only on $w.$ On the other hand, if we put
$\tilde{D}=\{z\in\mathbb{D}:\:|1-rz|\le\frac{3}{2}|1-r|\}$, then
\[
\int_{D}|f_{n}(z)|^{p}w(z)\mbox{d}A(z)\ge\frac{1}{(3/2)^{pn}(1-r)^{pn}}\int_{\tilde{D}}w(z)\mbox{d}A(z).
\]
 Since $r$ (thus $D$ and $\tilde{D}$) are fixed we see that
\[
\frac{1}{2^{pn}(1-r)^{pn}}=o\bigg(\frac{1}{(3/2)^{pn}(1-r)^{pn}}\int_{\tilde{D}}w(z)\mbox{d}A(z)\bigg),\qquad n\to\infty.
\]
 Thus,
\[
\frac{\|f'_{n}\|_{L_{a}^{p}(w)}^{p}}{\|f_{n}\|_{L_{a}^{p}(w)}^{p}}
\sim\int_{D}|f'_{n}(z)|^{p}w(z)\mbox{d}A(z)\Big/\int_{D}|f_{n}(z)|^{p}w(z)\mbox{d}A(z).
\]
Obviously,
\begin{align*}
\int_{D}|f'_{n}(z)|^{p}w(z)\mbox{d}A(z)= & \int_{D}\frac{n^{p}r^{p}}{|1-rz|^{pn+p}}w(z)\mbox{d}A(z)\\
\geq & \frac{n^{p}r^{p}}{2^{p}(1-r)^{p}}\int_{D}\frac{1}{|1-rz|^{pn}}w(z)\mbox{d}A(z)\\
= & \frac{n^{p}r^{p}}{2^{p}(1-r)^{p}}\int_{D}|f_{n}(z)|^{p}w(z)\mbox{d}A(z).
\end{align*}
Thus,
\[
\liminf_{n\to\infty}\frac{\|f'_{n}\|_{L_{a}^{p}(w)}}
{n\|f_{n}\|_{L_{a}^{p}(w)}}\ge\frac{r}{2(1-r)}.
\]
\end{proof}

For the proof of Theorem \ref{ourth2} we will need two lemmas.


\begin{lemma}
\label{lemma1}
Let $r\in[0,\,1)$ and $t\geq0.$ We set
\[
  I(t,\, r)=\int_{\mathbb{T}}\left|1-r\xi\right|^{-t}{\rm d}m(\xi)\quad
and \quad \varphi_{r}(t)=\int_{\mathbb{T}}\left|1+r\xi\right|^{t}{\rm d}m(\xi).
\]
Then,
\[
  I(t,\, r)=\frac{1}{\left(1-r^{2}\right)^{t-1}}\varphi_{r}(t-2)
\]
for every $t\geq2,$ and  $t\mapsto\varphi_{r}(t)$ is an increasing
function on $[0,\,+\infty)$ for every $r\in[0,\,1).$ Moreover, both
\[
  r\mapsto\varphi_{r}(t-2)\quad and\quad r\mapsto I(t,\, r),
\]
are increasing on $[0,\,1)$, for all $t\geq0.
$\end{lemma}

\begin{proof}
Indeed, supposing that $t\geq2,$ we can write
\[
I(t,\, r)=\frac{1}{1-r^{2}}\int_{\mathbb{T}}\left|b_{r}'(\xi)\right|
\frac{1}{\left|1-r\xi\right|^{t-2}}{\rm d}m(\xi),
\]
where $b_{r}(z)=\frac{r-z}{1-rz})$. Using the fact that $b_r\circ b_r(z) =z$
and changing the variable in the above integral we get
\begin{eqnarray*}
I(t,\, r) & = & \frac{1}{1-r^{2}}
\int_{\mathbb{T}}\left|b_{r}'(\xi)\right|\frac{1}{\left|1-rb_{r}\circ b_{r}(\xi)\right|^{t-2}}{\rm d}m(\xi)\\
 & = & \frac{1}{1-r^{2}}\int_{\mathbb{T}}\frac{1}{\left|1-rb_{r}(\xi)\right|^{t-2}}{\rm d}m(\xi)\\
 & = & \frac{1}{\left(1-r^{2}\right)^{t-1}}\varphi_{r}(t-2),
\end{eqnarray*}
since $1-rb_{r}(z)=\frac{1-rz-r(r-z)}{1-rz}=\frac{1-r^{2}}{1-rz}.$
Now,
\[
\varphi_{r}(t)=\int_{0}^{2\pi}\exp\left(\frac{t}{2}\ln\left(1+r^{2}-2r\cos s\right)\right){\rm d}s,
\]
\[
\varphi_{r}'(t)=\frac{1}{4}\int_{0}^{2\pi}\ln\left(1+r^{2}+2r\cos s\right)\exp\left(\frac{t}{2}\ln\left(1+r^{2}+2r\cos s\right)\right){\rm d}s,
\]
and
\[
\varphi_{r}''(t)=\frac{1}{4}\int_{0}^{2\pi}\left[\ln\left(1+r^{2}-2r\cos s\right)\right]^{2}\exp\left(\frac{t}{2}\ln\left(1+r^{2}-2r\cos s\right)\right){\rm d}s\geq0,
\]
for every $t\geq0,$ $r\in[0,\,1).$ Thus, $\varphi_{r}$ is a convex
fonction on $[0,\,\infty)$ and $\varphi_{r}'$ is increasing on $[0,\,\infty)$
for all $r\in[0,\,1).$ Moreover,
\[
\varphi_{r}'(0)=\frac{1}{4}\int_{0}^{2\pi}\ln\left(1+r^{2}-2r\cos s\right){\rm d}s =0.
\]
Thus,
\[
\varphi_{r}'(t)\geq\varphi_{r}'(0)=0,\:\forall t\in[0,\,\infty),\, r\in[0,\,1),
\]
and so $\varphi_{r}$ is increasing on $[0,\,\infty)$. The fact that
\[
r\mapsto I(t,\, r),
\]
is increasing on $[0,\,1)$ for all $t\geq0$ is obvious since
\[
I(t,\, r)=\left\Vert \frac{1}{(1-rz)^{t/2}}\right\Vert _{H^{2}}^{2}
=\sum_{k\geq0}a_{k}^2 (t) r^{2k},
\]
where $a_{k}(t)$ is the $k$th Taylor coefficient of $(1-z)^{-t/2}.$
The same reasoning gives that $r\mapsto\varphi_{r}(t)$ is increasing
on $[0,\,1)$.
\end{proof}


\begin{lemma}
\label{lemma2}
If for some $r_{0}\in[0,\,1)$ and $\gamma>0$ the function
$\frac{w(\rho)}{\left(1-\rho^{2}\right)^{\gamma}}$
is increasing on $[r_{0},\,1)$, then
\[
\int_{r}^{1}\rho w(\rho)I(t,\, r\rho){\rm d}\rho
\asymp\int_{r_{0}}^{1}\rho w(\rho)I(t,\, r\rho){\rm d}\rho,
\]
for all $t$ such that $t\ge \gamma + 3$ and for all $r\geq r_{0}$,
with constants independent on $t$.
\end{lemma}

\begin{proof}
Clearly,
\[
\int_{r_{0}}^{1}\rho w(\rho)I(t,\, r\rho){\rm d}\rho
\geq\int_{r}^{1}\rho w(\rho)I(t,\, r\rho){\rm d}\rho,
\qquad r\in [r_0, 1).
\]
Moreover,
\[
\int_{r_{0}}^{1}\rho w(\rho)I(t,\, r\rho){\rm d}\rho=\int_{r}^{1}\rho w(\rho)I(t,\, r\rho){\rm d}\rho+\int_{r_{0}}^{r}\rho w(\rho)I(t,\, r\rho){\rm d}\rho,
\]
 and applying Lemma \ref{lemma1},
\begin{eqnarray*}
\int_{r_{0}}^{r}\rho w(\rho)I(t,\, r\rho){\rm d}\rho & = & \int_{r_{0}}^{r}\frac{\rho w(\rho)}{\left(1-\rho{}^{2}\right)^{\gamma}}\frac{\left(1-\rho^{2}\right)^{\gamma}}{\left(1-(r\rho)^{2}\right)^{t-1}}\varphi_{r\rho}(t){\rm d}\rho\\
 & \leq & \frac{w(r)}{\left(1-r^{2}\right)^{\gamma}}\int_{r_{0}}^{r}\frac{\rho\left(1-\rho^{2}\right)^{\gamma}}{\left(1-(r\rho)^{2}\right)^{t-1}}\varphi_{r\rho}(t){\rm d}\rho\\
 & \leq & \frac{w(r)}{\left(1-r^{2}\right)^{\gamma}}\varphi_{r^{2}}(t)\int_{r_{0}}^{1}\frac{\rho\left(1-\rho^{2}\right)^{\gamma}}{\left(1-(r\rho)^{2}\right)^{t-1}}{\rm d}\rho,
\end{eqnarray*}
because $u\mapsto\varphi_{u}(t)$ is increasing for all $t>0$. For
the same reason,
\begin{eqnarray*}
 \int_{r}^{1}\rho w(\rho)\frac{1}{\left(1-(r\rho)^{2}\right)^{t-1}}
 \varphi_{r\rho}(t){\rm d}\rho & =
 & \int_{r}^{1}\frac{w(\rho)}{\left(1-\rho{}^{2}\right)^{\gamma}}
 \frac{\rho\left(1-\rho^{2}\right)^{\gamma}}{(1-(r\rho)^{2})^{t-1}}
 \varphi_{r\rho}(t){\rm d}\rho \\
 & \geq & \frac{w(r)}{\left(1-r^{2}\right)^{\gamma}}
 \varphi_{r^{2}}(t)\int_{r}^{1}\frac{\rho\left(1-\rho^{2}\right)^{\gamma}}{\left(1-(r\rho)^{2}\right)^{t-1}}{\rm d}\rho.
\end{eqnarray*}
Now note that
\[
\int_{r_0}^{r}\frac{\rho\left(1-\rho^{2}\right)^{\gamma}}
{(1-(r\rho)^{2})^{t-1}}{\rm d}\rho
\lesssim  \int_{r}^{1}
\frac{\rho\left(1-\rho^{2}\right)^{\gamma}}{(1-(r\rho)^{2})^{t-1}}{\rm d}\rho,
\qquad r\in [r_0, 1),
\]
with constants independent on $t \ge \gamma+3$.
Indeed, this estimate holds for $t=\gamma+3$, and, hence,
by monotonicity of the function
$\rho \mapsto (1-(\rho r)^2)^{-1}$, for all $ t\ge \gamma+3$.

Thus, using Lemma \ref{lemma1} and the fact that the function
$(1-\rho)^{-\gamma}w(\rho)$
is increasing on $[r_0,1)$, we obtain
\begin{eqnarray*}
\int_{r_{0}}^{r}\rho w(\rho)I(t,\, r\rho){\rm d}\rho &
\leq & \frac{w(r)}{\left(1-r^{2}\right)^{\gamma}}\varphi_{r^{2}}(t-2)
\int_{r_{0}}^{r}\frac{\rho\left(1-\rho^{2}\right)^{\gamma}}
{\left(1-(r\rho)^{2}\right)^{t-1}}{\rm d}\rho
\\
& \leq & \kappa_{1}\frac{w(r)}{\left(1-r^{2}\right)^{\gamma}}
\varphi_{r^{2}}(t-2)\int_{r}^{1}\frac{\rho\left(1-\rho^{2}\right)^{\gamma}}{\left(1-(r\rho)^{2}\right)^{t-1}}{\rm d}\rho\\
 & \leq & \kappa_{2}\int_{r}^{1}\rho w(\rho)\frac{1}{\left(1-(r\rho)^{2}\right)^{t-1}}\varphi_{r\rho}(t){\rm d}\rho,
\end{eqnarray*}
(where $\kappa_{1},\,\kappa_{2}$ are positive constants which do
not depend on $t$), which completes the proof.
\end{proof}

\begin{proof}[Proof of Theorem \ref{ourth2}]
We need to prove only the lower bound, the upper bound is already
proved in Theorem \ref{ourth1}.
Let us prove the minoration with the test function
$f(z)=\frac{1}{(1-rz)^{n}}$. Using \eqref{8} with $\alpha=r_{0}$,
we need to show that
\[
\begin{aligned}
\frac{\|f'\|^p_{L^p_a(w)}}{n^p r^p} & = \int_{r_{0}}^{1}\rho w(\rho)I(pn+p,\,
r\rho){\rm d}\rho \\
& \ge\frac{C}{(1-r)^{p}}
\int_{r_{0}}^{1}\rho w(\rho)I(pn,\, r\rho){\rm d}\rho =
\frac{C}{(1-r)^{p}}\|f\|^p_{L^p_a(w)}.
\end{aligned}
\]
Since $r\in [r_0, 1)$ and $n\ge \frac{\gamma+3}{p}$,
by Lemma \ref{lemma2} applied with $t=pn+p$ and $t=pn$ this means that
\[
\int_{r}^{1}\rho w(\rho)I(pn+p,\, r\rho){\rm d}\rho
\ge\frac{C}{(1-r)^{p}}\int_{r}^{1}\rho w(\rho)I(pn,\, r\rho){\rm d}\rho.
\]
By Lemma \ref{lemma1}, this is equivalent to the estimate

\begin{multline*}
\int_{r}^{1}\rho w(\rho)\frac{\varphi_{r\rho}(pn+p-2)}
{\left(1-(r\rho)^{2}\right)^{pn+p-1}}
{\rm d}\rho
\\
\geq\frac{C}{(1-r)^{p}}\int_{r}^{1}\rho w(\rho)
\frac{\varphi_{r\rho}(pn-2)}{\left(1-(r\rho)^{2}\right)^{pn-1}}{\rm d}\rho.
\end{multline*}
The last statement is obvious since

\begin{multline*}
\int_{r}^{1}\rho w(\rho)\frac{\varphi_{r\rho}(pn+p-2)}
{\left(1-(r\rho)^{2}\right)^{pn+p-1}}{\rm d}\rho   \\
\geq \frac{1}{(1-r^{2})^{p}}
\int_{r}^{1}\rho w(\rho)\frac{\varphi_{r\rho}(pn+p-2)}
{\left(1-(r\rho)^{2}\right)^{pn-1}}{\rm d}\rho \\
\geq \frac{1}{(1-r^{2})^{p}}\int_{r}^{1}
\rho w(\rho)\frac{\varphi_{r\rho}(pn-2)}{\left(1-(r\rho)^{2}\right)^{pn-1}}{\rm d}\rho,
\end{multline*}
where the last inequality is due to the fact that $t\mapsto\varphi_{u}(t)$
is increasing for all $0\leq u<1$. \end{proof}


\section{The case of super-polynomially
decreasing weights. Proof of Theorem \ref{ourth3}:}

For the proof of Theorem \ref{ourth3} we will need a definition from
the theory of model subspaces of the Hardy space.
For a finite subset $\sigma$ of $\mathbb{D}$
with ${\rm card}\,\sigma=n$, consider the finite Blaschke product
\[
B_{\sigma}=\prod_{\lambda\in\sigma}b_{\lambda},
\]
where $b_{\lambda}(z)=\frac{\lambda-z}{1-\overline{\lambda}z}$,
$\lambda\in\mathbb{D}$.
Define the model space $K_{B_{\sigma}}$ by
$$
  K_{B_{\sigma}} = \left(B_{\sigma}H^{2}\right)^{\perp}=H^{2}
  \ominus B_{\sigma}H^{2}.
$$
Consider the family $\left(e_{k}\right)_{1\leq k\leq n}$ in $K_{B_{\sigma}}$
(known as {\it Malmquist basis}, see \cite[p. 117]{Nik}),
$$
  e_{1}(z)=\frac{(1-|\lambda_1|)^{1/2}}
  {1-\overline{\lambda_1} z}\quad \mbox{and} \quad
  e_{k}(z)=\bigg({\displaystyle \prod_{j=1}^{k-1}}b_{\lambda_{j}}(z)\bigg)
  \frac{(1-|\lambda_k|)^{1/2}}{1 - \overline{\lambda_k} z}, \quad k\in [2,\, n],
$$
The family $\left(e_{k}\right)_{1\leq k\leq n}$
associated with $\sigma$ is an orthonormal basis of the $n$-dimensional
space $K_{B_{\sigma}}.$

In what follows we denote by $L_{a}^{p}(w,\, s\mathbb{D})$
and by $H^{p}(s\mathbb{D})$, $s>0$, the weighted Bergman space
and the Hardy space in the disc $s\mathbb{D} = \{z: |z|<s\}$,
respectively. If $w\equiv 1,$ we write simply
$L_{a}^{p}(s\mathbb{D})$ and we write $L^p_a$ if $s=1$.


\begin{lemma}
\label{lemma3}
Let $n\geq1,$ $r,\, s\in[0,\,1)$
and $p\in[1,\,+\infty].$ We set
\[
   M_{p,\, s}(n,\, r) = \sup\Big\{ |f(\xi)|:\:\xi\in\mathbb{D},\,
   f\in\mathcal{R}_{n,\, r},\,\left\Vert f\right\Vert _{L_{a}^{p}
   \left(s\mathbb{D}\right)} \leq 1 \Big\}.
\]
Then
\begin{equation}
\label{11}
  M_{p,\,\frac{2}{3}}(n,\, r)\le d\frac{c^{n}}{(1-r)^{n+b}},
\end{equation}
where $d>0,\, b>0$, $c>1$ are some absolute positive constants
$($may be, depending on $p$$)$.
\end{lemma}

\begin{remark} Lemma \ref{lemma3}
is valid not only for $s=\frac{2}{3}$, but for every $s\in(0,\,1),$
with constants $d>0,\, b>0$, $c>1$ depending both on $s$ and $p$.
\end{remark}

\begin{proof} For every $f\in\mathcal{R}_{n,\, r}$
and $\xi \in \mathbb{D}$, we have
\[
\left|f\left(\frac{1}{2}\xi\right)\right|=
\left|f_{\frac{2}{3}}\left(\frac{3}{4}\xi\right)\right|=
\left|\int_{\mathbb{D}} f_{\frac{2}{3}} \left(u\right)
\left(\overline{k_{\frac{3}{4}\xi}(u)}\right)^{2}{\rm d}A(u)\right|,
\]
where $k_{\lambda}(z)=\frac{1}{1-\overline{\lambda}z}$
is the standard Cauchy kernel associated with $\lambda\in\mathbb{D}$,
and $A$ is the normalized area measure on $\mathbb{D}$.
Applying H\"older's inequality we obtain
\[
   \left|f\left(\frac{1}{2}\xi\right)\right|
   \leq  \left\Vert f_{\frac{2}{3}}\right\Vert _{L_{a}^{p}}
   \left\Vert \left(k_{\frac{3}{4}\xi}\right)^{2}\right\Vert _{L_{a}^{p'}}
   = \left(\frac{3}{2}\right)^{1/p}
   \left\Vert f\right\Vert _{L_{a}^{p}
   \left(\frac{2}{3}\mathbb{D}\right)}\left\Vert
   \left(k_{\frac{3}{4}\xi}\right)^{2}\right\Vert _{L_{a}^{p'}},
   \quad \xi \in \mathbb{T},
\]
where $p'$ is such that $\frac{1}{p}+\frac{1}{p'}=1$.
Now, note that
\[
 \left\Vert \left(k_{\frac{3}{4}\xi}\right)^{2}
\right\Vert _{L_{a}^{p'}}\leq\left\Vert
\left(k_{\frac{3}{4}\xi}\right)^{2}\right\Vert_{H^{\infty}}
= \left(\frac{1}{1-\frac{3}{4}}\right)^{2} = 16.
\]
Finally, supposing $\left\Vert f\right\Vert _{L_{a}^{p}
\left(\frac{2}{3}\mathbb{D}\right)}\leq1$, we obtain
\[
\left\Vert f\right\Vert _{L_{a}^{p}\left(\frac{1}{2}\mathbb{D}\right)}\leq\left\Vert f\right\Vert _{H^{\infty}\left(\frac{1}{2}\mathbb{D}\right)}
\leq 16\left(\frac{3}{2}\right)^{1/p} \leq 24,
\]
which gives

\begin{equation}
\label{11-1}
M_{p,\,\frac{2}{3}}(n,\, r)\leq24M_{2,\,\frac{1}{2}}(n,\, r).
\end{equation}
\par

It remains to obtain a suitable upper bound for $M_{2,\,\frac{1}{2}}(n,\, r).$
Let us prove that

\begin{equation}
\label{12}
M_{2,\,\frac{1}{2}}(n,\, r) \leq 2
\sqrt{n}\left(\frac{2}{1-r}\right)^{n+\frac{1}{2}}.
\end{equation}
For every $f\in\mathcal{R}_{n,\, r}$, we have  $f_{\frac{1}{2}}\in
\mathcal{R}_{n,\,\frac{1}{2}r}\subset\mathcal{R}_{n,\, r}$.
If $\{1/\overline{\lambda}_1, \dots, 1/\overline{\lambda}_{n}\}$
is the set of the poles of $f$
(thus, $|\lambda_j| \leq r$, $j=1, \dots, n$),
then $f\in K_{B_{\sigma}}$ with
$\sigma=\{ \lambda_1, \dots, \lambda_n\} \subset r\mathbb{D}$),
whereas the set
$\{2/\overline{\lambda}_1,\dots , 2/\overline{\lambda}_n\}$
is the set of the poles of the function $f_{\frac{1}{2}}$
and $f_{\frac{1}{2}}\in K_{B_{\sigma'}}$
with $\sigma'=\left\{ \frac{1}{2}\lambda_1,\dots ,\frac{1}{2}\lambda_n
\right\} \subset\frac{r}{2}\mathbb{D}$.
Hence, there exist $\alpha_{1}, \dots, \alpha_n \in\mathbb{C}$
such that
\begin{equation}
\label{13}
f_{\frac{1}{2}} = \sum_{k=1}^{n} \alpha_{k} e_{k},
\end{equation}
on $\mathbb{D}$, where $(e_k)_{k=1}^n$ is the Malmquist basis
associated with the set $\sigma'$. Since both $f_{\frac{1}{2}}$
and $\sum_{k=1}^{n}\alpha_{k} e_{k}$ are meromorphic in $\mathbb{C}$
the equality \eqref{13} is in fact valid everywhere in $\mathbb{C}$.
Thus,
\[
f\left(\xi\right)=\sum_{k=1}^{n}\alpha_{k}
\left({\displaystyle \prod_{j=1}^{k-1}}
\frac{\frac{\lambda_{j}}{2}-2\xi}{1-\overline{\lambda_{j}}\xi}\right)
\frac{\left(1-\frac{1}{4}\left|\lambda_{k}\right|^{2}\right)^{1/2}}
{1-\overline{\lambda_{j}}\xi}, \qquad \xi\in\mathbb{D},
\]
and by the Cauchy--Schwarz inequality,
\begin{equation}
\label{12a}
\left|f\left(\xi\right)\right|\leq
\Big(\sum_{k=1}^{n} |\alpha_{k}|^2 \Big)^{1/2}
\Bigg(\sum_{k=1}^{n}\Bigg| \bigg({\displaystyle
\prod_{j=1}^{k-1}}\frac{\frac{\lambda_{j}}{2}-2\xi}
{1-\overline{\lambda}_j \xi} \bigg)
\frac{\big(1 -\frac{1}{4}|\lambda_{k}|^{2} \big)^{1/2}} {1-\overline{\lambda_{j}}\xi}
\Bigg|^{2}\Bigg)^{1/2}.
\end{equation}
for any $\xi\in\mathbb{D}$.  Now, if $\lambda\in
r\mathbb{D}$ and $\xi \in \mathbb{D}$,
\[
\frac{\frac{\lambda}{2}-2\xi}{1-\overline{\lambda}\xi}
=\frac{2\left(\frac{\lambda}{4}-\xi\right)}
{1-\frac{\overline{\lambda}}{4}\xi}
\frac{1-\frac{\overline{\lambda}}{4}\xi}{1-\overline{\lambda}\xi}=
2 b_{\frac{\lambda}{4}} (\xi)  \bigg( 1+\frac{3\overline{\lambda}}
{4(1-\overline{\lambda}\xi)}\bigg),
\]
which gives
\[
\bigg|\frac{\frac{\lambda}{2}-2\xi}{1-\overline{\lambda}\xi}\bigg|
\leq  2\left(1+\frac{3r}{4}\frac{1}{1-r}\right)
= \frac{4-r}{2(1-r)} \leq \frac{2}{1-r}.
\]
We get
\begin{equation}
\label{12b}
\begin{aligned}
 \sum_{k=1}^{n} \Bigg|\bigg({\displaystyle \prod_{j=1}^{k-1}}
 \frac{\frac{\lambda_{j}}{2}-2\xi}{1-\overline{\lambda_{j}}\xi}
 \bigg) \frac{\big(1-\frac{1}{4} |\lambda_{k}|^{2}\big)^{1/2}}
 {1-\overline{\lambda_{j}}\xi}\Bigg|^{2}
 & \leq \frac{1}{(1-r)^{2}}\sum_{k=1}^{n}2^{2(k-1)}
 \left(\frac{1}{1-r}\right)^{2(k-1)}\\
 & \leq \frac{1}{4}\left(\frac{2}{1-r}\right)^{2n+1}.
\end{aligned}
\end{equation}

\par
Now we first notice that
\[
  \Big(\sum_{k=1}^{n}\left|\alpha_{k}\right|^{2}\Big)^{1/2}=
\left\Vert f_{\frac{1}{2}}\right\Vert _{H^{2}}.
\]
For any function $\varphi(z) = \sum_{k\ge 0} \widehat{\varphi}(k) z^k$
in $H^2$, one has
\[
  \left\Vert \varphi\right\Vert _{H^{2}}^{2}=
  \sum_{k\geq0}\frac{\left|\widehat{\varphi}(k)\right|}
  {\sqrt{k+1}}\sqrt{k+1}\left|\widehat{\varphi}(k)\right|
  \leq  \left\Vert \varphi\right\Vert _{L_{a}^{2}}
  \left\Vert \varphi\right\Vert_{B_{2,\,2}^{1/2}},
\]
We now use the upper bound of  \cite[Theorem A, (4)]{Z6}:
for $\varphi \in \mathcal{R}_{n,\,\rho}$ one has
\begin{eqnarray*}
\|\varphi\|_{B_{2,\,2}^{1/2}}^{2} & = & \|\varphi'\|_{L_{a}^{2}}^{2}+\|\varphi\|_{H^{2}}^{2}\\
 & \leq & (2+r)\frac{n}{1-r}\|\varphi\|_{H^{2}}^{2}+\|\varphi\|_{H^{2}}^{2}\\
 & \leq & \frac{4n}{1-r}\|\varphi\|_{H^{2}}^{2},\end{eqnarray*}
which gives
\[
\left\Vert \varphi\right\Vert _{H^{2}}\leq2\sqrt{\frac{n}{1-\rho}}\left\Vert \varphi\right\Vert _{L_{a}^{2}}.
\]
In particular, with $\varphi=f_{\frac{1}{2}}$
we get $\varphi\in\mathcal{R}_{n,\,\frac{1}{2}r}$ and
\begin{equation}
\label{12c}
\left\Vert f_{\frac{1}{2}}\right\Vert _{H^{2}} \leq
2\sqrt{\frac{n}{1-r/2}}\left\Vert f_{\frac{1}{2}}\right\Vert _{L_{a}^{2}}
\leq 2 \sqrt{2n} \left\Vert f_{\frac{1}{2}}\right\Vert _{L_{a}^{2}}.
\end{equation}
We conclude from (\ref{12a}), (\ref{12b}) and (\ref{12c})
that for any $\xi\in\mathbb{D}$,
\[
\left|f\left(\xi\right)\right|\leq
\left\Vert f_{\frac{1}{2}}\right\Vert _{H^{2}}\left(\frac{1}{4}
\frac{2^{2n+1}}{(1-r)^{2n+1}}\right)^{\frac{1}{2}}\leq\frac{1}{2}
\left(\frac{2}{1-r}\right)^{n+\frac{1}{2}}2 \sqrt{2n}
\left\Vert f_{\frac{1}{2}}\right\Vert _{L_{a}^{2}},
\]
that is,
\[
\left|f\left(\xi\right)\right|
\leq \sqrt{2n}\left(\frac{2}{1-r}\right)^{n+\frac{1}{2}}
\left\Vert f\right\Vert _{L_{a}^{2}\left(\frac{1}{2}\mathbb{D}\right)},
\qquad \xi\in\mathbb{D}.
\]
Taking the supremum over $\xi\in\mathbb{D}$ and
$f\in\mathcal{R}_{n,\, r}$ we obtain \eqref{12}.

Combining \eqref{11-1} and \eqref{12} and choosing $d=48$,
$b=\frac{1}{2}$ and $c>2$ such that $2^n \sqrt{n} \leq c^{n}$ for any $n\geq1$,
we complete the proof and obtain \eqref{11}.
\end{proof}


\begin{proof}[Proof of Theorem \ref{ourth3}]
Take $r\in(0,\,1)$ and $R\in(0,\, r)$ and let us
represent the norm $\|f'\|_{L_{a}^{p}(w)}^{p}$
of a function $f \in \mathcal{R}_{n,\,r}$ as $I_{1}+I_{2}$,
\[
I_{1}=\int_{0}^{R}\|(f_{\rho})'\|_{p}^{p}w(\rho)\mbox{d}\rho,\qquad I_{2}=\int_{R}^{1}\|(f_{\rho})'\|_{p}^{p}w(\rho)\mbox{d}\rho.
\]
Here and everywhere below in this proof, $C_{i},\, i=1, \dots, 5,$ are
positive constants, depending, may be, only on $p$ and $w$
(but not on $n$ and $r$).
By (\ref{hardy}), we have for the first integral
\[
I_{1}\le C_{1}\Big(\frac{n}{1-R}\Big)^{p}\int_{0}^{R}
\|f_{\rho}\|_{p}^{p}w(\rho)\mbox{d}\rho\le C_{2}
\Big(\frac{n}{1-R}\Big)^{p}\|f\|_{L_{a}^{p}(w)}^{p}.
\]
Note that $f_\rho \in \mathcal{R}_{n,\,\rho r} \subset \mathcal{R}_{n,\,r}$,
and, thus,
$\|f_\rho\|_\infty \le M_{p,\,\frac{2}{3}}(n,\, r)
\|f_\rho\|_{L^p_a(\frac{2}{3}\mathbb{D})}$.
Applying (\ref{hardy}) once again together with an obvious inequality
$\|f_\rho\|_p \le \|f_\rho\|_\infty$,
we get
\[
\begin{aligned}
I_{2} & \le C_{3}\Big(\frac{n}{1-r}\Big)^{p}\int_{R}^{1}
\|f_{\rho}\|_{\infty}^{p}w(\rho)\mbox{d}\rho \\
      & \le C_{3}\|f\|_{L_{a}^{p}\left(\frac{2}{3}\mathbb{D}\right)}^{p}
      \Big(\frac{n}{1-r}\Big)^{p}\int_{R}^{1}M^p_{p,\,2/3}(n,\, r)
      w(\rho)\mbox{d}\rho \\
      & \le C_{3}\|f\|_{L_{a}^{p}\left(\frac{2}{3}\mathbb{D}\right)}^{p}
      \Big(\frac{n}{1-r}\Big)^{p}\frac{c^{pn}}{(1-r)^{pn+pb}}w(R),
\end{aligned}
\]
where the last inequality follows from Lemma \ref{lemma3}.
Note that \[\|f\|_{L_{a}^{p}\left(\frac{2}{3}\mathbb{D}\right)}^{p}\le(w(2/3))^{-1}
\|f\|_{L_{a}^{p}(w)}^{p}.\] Hence,
\[
  I_{2}\le C_{4}\Big(\frac{n}{1-r}\Big)^{p}\frac{c^{pn}}{(1-r)^{pn+pb}}w(R)\|f\|_{L_{a}^{p}(w)}^{p}.
\]
Now, choose a positive increasing sequence
$\left(\gamma_{n}\right)_{n\in\mathbb{N}}$
such that $n=o\left(\gamma_{n}\right)$, as $n \to +\infty$.
For any $n$ we fix $r_n^\circ$ such that the function
$w(r)(1-r)^{-\gamma_{n}}$ decreases on $[r_n^\circ,\,1)$.
Now for a fixed $n$ take $r,\, R$ so that $r_n^\circ < R < r<1$ and
\[
  1-R=(1-r)^{1/2}, \qquad 1-r_n^\circ = (1-r)^{1/4}.
\]
We have
\[
  w(R)\le w(r_n^\circ) \frac{(1-R)^{\gamma_{n}}}
  {(1-r_n^\circ)^{\gamma_{n}}} = w(r_n^\circ) (1-r)^{\gamma_{n}/4}.
\]
Hence, using the fact that $w$ is bounded on $[r_1^\circ, 1)$, we obtain
\[
  I_{2}\le C_{4}\Big(\frac{n}{1-r}\Big)^{p}\|f\|_{L_{a}^{p}(w)}^{p}\cdot c^{pn}
  \frac{(1-r)^{\gamma_{n}/4}}{(1-r)^{pn+pb}}.
\]
Let us show that for sufficiently large $n$,
$$
  c^{pn}\frac{(1-r)^{\gamma_{n}/4}}{(1-r)^{pn+pb}} \to 0,\qquad r\to1-.
$$
Indeed, choosing $r$ so that $c<(1-r)^{-1}$, we get
\[
c^{pn}\frac{(1-r)^{\gamma_{n}/4}}{(1-r)^{pn+pb}}
\le(1-r)^{\frac{\gamma_{n}}{4}-2pn-pb}\to0,\qquad r\to1-,
\]
since $n = o(\gamma_{n})$, $n\to\infty$. Hence,
there exists a sequence $(r_{n})$, $r_n \to 1-$, such that
\[
\frac{I_{2}^{1/p}}{n\|f\|_{L_{a}^{p}(w)}}=o\bigg(\frac{1}{1-r_{n}}\bigg),
\qquad n\to \infty.
\]
The corresponding estimate for $I_{1}$ is obvious since
$1-R_{n}=(1-r_{n})^{1/2}$.
\end{proof}

\bibliographystyle{amsplain}

\providecommand{\bysame}{\leavevmode\hbox to3em{\hrulefill}\thinspace}
\providecommand{\MR}{\relax\ifhmode\unskip\space\fi MR }
\providecommand{\MRhref}[2]{%
  \href{http://www.ams.org/mathscinet-getitem?mr=#1}{#2}
}
\providecommand{\href}[2]{#2}
\begin{thebibliography}{}

\end{thebibliography}


\begin{thebibliography}{References}


\bibitem{bar1} A.D. Baranov, \textit{Bernstein-type
inequalities for shift-coinvariant subspaces and their
applications to Carleson embeddings},
J. Funct. Anal. {\bf 223} (2005), no. 1, 116--146.

\bibitem{bar2} A.D. Baranov, \textit{Embeddings of model subspaces
of the Hardy space: compactness and Schatten--von Neumann ideals},
{\it Izvestia RAN Ser. Matem.} {\bf 73} (2009), no. 6, 3--28;
English transl.: {\it Izv. Math.} {\bf 73} (2009), no. 6, 1077--1100.

\bibitem{BoEr} P. Borwein and T. Erd\'elyi, \textit{Polynomials and
Polynomial Inequalities}, Springer, New York, 1995.

\bibitem{Da} V.I. Danchenko, \textit{An integral estimate for the
derivative of a rational function}, Izv. Akad. Nauk SSSR Ser. Mat.
{\bf 43} (1979), no. 2, 277--293;
English transl. Math. USSR Izv., {\bf 14} (1980), no. 2, 257--273.

\bibitem{Dol} E.P. Dolzhenko, \textit{Rational approximations and boundary properties of analytic functions},
  Mat. Sb. (N.S.), 69(111):4 (1966), 497–-524.


\bibitem{Dy2} K.M. Dyakonov, \textit{Kernels of Toeplitz operators,
smooth functions, and Bernstein-type inequalities}, Zap. Nauchn.
Semin. S. Peterburg. Otdel. Mat. Inst. Steklov. (POMI) {\bf 201} (992),
5--21; English transl.: J. Math. Sci. {\bf 78} (1996), 131--141.

\bibitem{Dy1} K.M. Dyakonov, \textit{Differentiation in star-invariant
subspaces I. Boundedness and compactness,} J. Funct. Anal. {\bf 192 } (2002),
364--386.

\bibitem{ind} K.M. Dyakonov, \textit{Smooth functions in the range
of a Hankel operator}, Indiana Univ. Math. J. {\bf 43} (1994),
805--838.

\bibitem{Erd} T. Erd\'elyi, \textit{George Lorentz and inequalities
in approximation}, Algebra i Analiz  {\bf 21} (2009), no. 3, 1--57.

\bibitem{Go} A.A. Gonchar, \textit{Degree of approximation by rational
fractions and properties of functions}, Proc. Internat. Congr. Math.
(Moscow, 1966), Mir, Moscow, 1968, 329--356; English transl.:
Amer. Math. Soc. Transl. (2) {\bf 91} (1970).

\bibitem{HedKorZhu} H. Hedenmalm, B. Korenblum, and K. Zhu,
\textit{Theory of Bergman spaces}, Graduate Texts in Mathematics,
{\bf 199}, Springer-Verlag, New York, 2000.

\bibitem{LeTr} R.J. Leveque, L.N. Trefethen, \textit{On the resolvent
condition in the Kreiss matrix Theorem}, BIT {\bf 24} (1984), 584--591.

\bibitem{Nik} N. Nikolski, \textit{Treatise on the Shift Operator},
Springer-Verlag, Berlin, 1986.

\bibitem{Pek} A.A. Pekarskii, \textit{Inequalities of Bernstein
type for derivatives of rational functions, and inverse theorems of
rational approximation,} Mat. Sb. {\bf 124(166)} (1984), 571-588; English transl.:
Math. USSR-Sb. {\bf 52} (1985), no. 2, 557--574.

\bibitem{Pel} V.V. Peller, \textit{Hankel operators of class $\mathcal{S}_{p}$
and their applications (rational approximations, Gaussian processes,
the problem of majorizing operators)}, Mat. Sb. {\bf 113(155)} (1980), no. 4,
538--581; English transl.: Math. USSR Sb. {\bf 41} (1982), no. 4, 443--479.

\bibitem{Ru} V.N. Rusak, \textit{Rational Functions as Approximation
Apparatus}, Izdat. Beloruss. Gos. Univ., Minsk, 1979 (Russian).

\bibitem{Sp} M.N. Spijker, \textit{On a conjecture by LeVeque and
Trefethen related to the Kreiss matrix theorem,} BIT {\bf 31} (1991), 551--555.

\bibitem{Z5} R. Zarouf, \textit{Sharpening a result by E.B. Davies
and B. Simon}, C. R. Acad. Sci. Paris {\bf 347} (2009), no. 16, 939--942.

\bibitem{Z1} R. Zarouf, \textit{Analogs of the Kreiss resolvent condition
for power bounded matrices}, Actes des journ\'ees du GDR
"Analyse Fonctionnelle Harmonique et Applications", Metz, 2010.

\bibitem{Z2} R. Zarouf, \textit{Asymptotic sharpness of a Bernstein-type
inequality for rational functions in $H^{2}$,} Algebra i Analiz
{\bf 23} (2011), no. 2, 147--161.

\bibitem{Z6} R. Zarouf, \textit{Application of a Bernstein type inequality
to rational interpolation in the Dirichlet space},
Zap. Nauchn. Semin. S. Peterburg. Otdel. Mat. Inst. Steklov. (POMI)

39:101--112, (2011).

\bibitem{Z4} R. Zarouf, \textit{Effective $H^{\infty}$ interpolation
constrained by Hardy and Bergman weighted norms,} Ann. Funct. Anal.,
2(2):59--74 (2011).

\bibitem{Z3} R. Zarouf, \textit{Effective $H^{\infty}$ interpolation},
to appear in Houston J. Math.

\bibitem{Zhu} K. Zhu, \textit{Operator theory in function spaces},
Monographs and Textbooks in Pure and Applied Mathematics, 139. Marcel
Dekker, Inc., New York, 1990. \end{thebibliography}

\end{document}